\documentclass[a4paper,12pt]{article}
\usepackage{amsmath,amssymb,amsthm,mathrsfs}

\usepackage{pgfplots,pgfplotstable,filecontents}
\usepackage[labelfont=bf]{caption}
\usepackage{subcaption}

\usepackage{graphicx}
\usepackage{epstopdf}

\usepackage[left=1in,right=1in,top=1in,bottom=1in,footskip=.25in]{geometry}
\usepackage{fancyhdr}
\newtheorem{theorem}{Theorem}

\newtheorem{remark}{Remark}

\def\eqd{\,{\buildrel d \over =}\,}

\newcommand{\vnorm}[1]{\left\|#1\right\|}

\newcommand{\eop}{\hfill$\Box$\medskip}

\newcommand{\pu}[1]{\textcolor{black}{#1}}

\begin{document}

\title{Approximating the equilibrium quantity traded and welfare in large markets}

\author{Ellen V. Muir\thanks{The corresponding author. School of Mathematics and Statistics, The University of Melbourne, Parkville VIC 3010, Australia; e-mail: evmuir@unimelb.edu.au.} \  and Konstantin Borovkov\thanks{School of Mathematics and Statistics, The University of Melbourne, Parkville VIC 3010, Australia; e-mail: borovkov@unimelb.edu.au.}}

\date{}

\maketitle

\begin{abstract} \pu{We consider the efficient outcome of a canonical economic market model involving buyers and sellers with independent and identically distributed random valuations and costs, respectively. When the number of buyers and sellers is large, we show that the joint distribution of the equilibrium quantity traded and welfare is asymptotically normal. Moreover, we bound the approximation rate. The proof proceeds by constructing, on a common probability space, a representation consisting of two independent empirical quantile processes, which in large markets can be approximated by independent Brownian bridges. The distribution of interest can then be approximated by that of a functional of a Gaussian process. This methodology applies to a variety of mechanism design problems.}
\medskip


\medskip

\hphantom{2010 Mathematics Subject Classification:}
Secondary 91B26; 60F17; 60G15

\end{abstract}

\setlength{\baselineskip}{1.5\baselineskip}

\section{Introduction and main results}
\label{sec:1}


\pu{The problem of how to achieve efficient market outcomes has been of fundamental importance in economics, dating back to the pioneering work of L\'{e}on Walras. A significant and sometimes insurmountable barrier to designing efficient markets is the information required to set prices that equate market demand and supply. In many practical applications this information is privately held by strategic buyers and sellers (agents) who will not freely reveal this information, given the impact this will have on market prices. For a market intermediary (such as a stock exchange) to elicit this information and implement the efficient outcome, agents with private information must be appropriately compensated. Thus, in many cases it is not possible to design a market which is efficient, correctly elicits private information from agents (formally, is incentive compatible and individually rational) and does not require an intermediary to subsidise trade (formally, is deficit-free). In the field of mechanism design, this result is known as the impossibility of efficient trade and was first emphasised by Myerson and Satterthwaite~\cite{ms83}.}

\pu{Mechanism design provides an important and influential approach to dealing with market design problems involving private information, with applications ranging from the allocation of government resources (including land, mining rights, radio spectrum licenses and university places) to kidney exchange programs and advertisement placement in Internet search engines (see Loertscher, Marx and Wilkening~\cite{lmw15} for a recent survey). In the Bayesian mechanism design literature the private information of a given buyer or seller is modelled as a random variable, whose realisation is known only to that buyer or seller. All other buyers and sellers and the market intermediary know only the distribution of this random variable. The role of the intermediary is to choose the market mechanism. A mechanism determines, as a function of the private information revealed by agents, the identity of trading agents and the prices they pay and receive. In light of the Myerson--Satterthwaite impossibility theorem, a prominent strand of mechanism design literature focuses on designing markets that are efficient, incentive compatible, individually rational and deficit-free asymptotically, as the number of buyers and sellers grows large (see, for example, McAfee~\cite{mcafee92} and Rustichini, Satterthwaite and Williams~\cite{rsw94}). This requires approximating mechanism outcomes (such as the equilibrium quantity traded and welfare) in large markets.}

\pu{In this paper we consider the canonical mechanism design market model, known as the \emph{independent private values model}. We devise a general methodology for approximating mechanism outcomes in large markets under the independent private values assumption. To do this, we restrict attention to approximating the equilibrium quantity traded and welfare under the efficient mechanism (which, of course, is not necessarily deficit-free). We show that, as both the number of buyers and number of sellers tend to infinity, the joint distribution of the equilibrium quantity traded and welfare is asymptotically normal, and give an upper bound for the approximation rate. This is accomplished by constructing, on a common probability space, an empirical quantile process representing market demand and an independent empirical quantile process representing supply together with two independent and appropriately weighted Brownian bridges approximating the above-mentioned quantile processes. The distribution of interest can then be approximated by that of a functional of a Gaussian process.}

\pu{Several papers analyse the performance of mechanisms in large markets with independent private values, including Gresik and Satterthwaite~\cite{gs89}, McAfee~\cite{mcafee92}, Satterthwaite and Williams~\cite{sw02} and Rustichini, Satterthwaite and Williams~\cite{rsw94}. However, this literature has previously focused on computing the rate at which mechanism outcomes converge to efficiency. To the best of our knowledge, this paper is the first to compute higher order distributional approximations (of the central limit theorem-type) to mechanism outcomes. One advantage of this approach is that it enables the direct comparison of different mechanisms using the parameters of the approximating normal distributions. Furthermore, our approach immediately generalises to any mechanism which can be appropriately represented in terms of transformed empirical quantile functions, which covers a large class of mechanisms studied in the Bayesian mechanism design literature, as is discussed in detail in Section 2. The common probability space method also allows us to compute the covariance of mechanism outcomes, which is problematic if a rate of convergence approach is adopted and is important in some settings (see, for example, p. 447 of McAfee~\cite{mcafee92}). Finally, unlike papers (such as Gresik and Satterthwaite~\cite{gs89}) which consider sequences of markets with a fixed ratio of buyers and sellers, we formulate more general convergence results that apply to nets of markets in which this ratio is not necessarily fixed.}

\pu{The remainder of this paper is structured as follows. We now introduce the independent private values model, describe the problem of interest and state our main results. In Section 2, we discuss extensions and applications of these results. Proofs are included in Section 3.}

We consider a market in which units of a homogeneous, indivisible good are traded among $N$ buyers $ \{i_1,\dots,i_N\}$ and $M$ sellers $ \{j_1,\dots,j_M\}$. Each buyer is interested in purchasing one unit of the good, and each seller  has the capacity to produce and sell one unit of the good. Buyers are willing to purchase the good at a price not exceeding their respective private reservation valuations $V_1,\dots,V_N$ that are assumed to be independent and identically distributed random variables with a common distribution function $F$ with support $[a,b]$. Sellers are willing to produce at a price that is not less than their respective private production costs $C_1,\dots,C_M,$ which are also assumed to be independent and identically distributed  random variables, with a distribution function $G$ with support $[c,d]$. We assume that buyers valuations and sellers costs are independent of each other and call the $(N+M)$-tuple
\[
\mathfrak{R}_{(N,M)}:= (V_1,\dots,V_N;C_1,\dots,C_M)
\]
a realization of the market $\mathfrak{M} := \langle N, M, F, G \rangle$.

\pu{This independent private values model has been studied extensively in the literature on mechanism design and auction theory, see, for example, Myerson~\cite{myerson81}, Milgrom and Webber~\cite{mw82}, Chatterjee and Samuelson~\cite{cs83}, Myerson and Satterthwaite~\cite{ms83}, Gresik and Satterthwaite~\cite{gs89}, Williams~\cite{williams99}, Baliga and Vohra~\cite{bv03}, Muir~\cite{muir13} and Loertscher and Marx~\cite{lm16}. The results of our analysis immediately generalise to settings in which sellers produce multiple units (and similarly, buyers demand multiple units), provided the cost of production for each unit is independently drawn from the distribution $G$. We may also relax the assumption of identical distributions. For example, we can suppose a fixed proportion of sellers draw their costs from some distribution $G_1$ and the rest draw their costs from some distribution $G_2$, provided we use the appropriate mixture distribution in our asymptotic analysis. The independent private values model is robust to the introduction of a small amount of dependence among valuations and costs when prices are bounded in magnitude (see  Kosmospoulou and Williams~\cite{kw98}). Thus, theory and results concerning this model are useful for many practical applications, while results concerning models with dependent types are considered fragile.}

Within a market, the welfare generated by trade is defined as the sum of trading buyer valuations less the sum of trading seller costs:
\[
\sum_{i \in \mathcal{N}} V_{i} - \sum_{j \in \mathcal{M}} C_j,
\]
where   $\mathcal{N} \subset \{i_1,\dots,i_N\}$ and $\mathcal{M}\subset\{j_1,\dots,j_N\}$ are, respectively,  the subsets of buyers and sellers who trade in the market. Feasibility requires $\left|\mathcal{N}\right| \leq \left|\mathcal{M}\right|$. Market welfare can be thought of as aggregating the gains from trade for all market participants.

A market is said to be \emph{efficient} if the level of market welfare is always maximised for given buyer valuations and seller costs. The quantity traded in an efficient market is called the efficient quantity. To compute it, we form a demand curve by ordering buyer valuations $V_i$ from highest to lowest: $V_{[1]} \geq \dots \geq V_{[N]}$, and we form a supply curve by ordering seller costs $C_j$ from lowest to highest: $C_{(1)} \leq \dots \leq C_{(M)}$.  The efficient quantity is then given by
\begin{align}
	\label{eq:kdef}
K
 := \underset{0 \leq k \leq N \wedge M }{\arg \max} \sum_{i =1}^{k} (V_{[i]} - C_{(i)})
 \equiv
  \left| \{ k \in \{1,\dots,N \wedge M\}:V_{[k]}-C_{(k)} \geq 0 \} \right|.
\end{align}
{\em In an efficient market, the value of welfare is given by}
\begin{align}
	\label{eq:welfQ}
W := \sum_{i=1}^{K} (V_{[i]} - C_{(i)}).
\end{align}
Figure~\ref{sd} provides an illustration for the quantities~$K$ and~$W$. \pu{Note that it is not necessary to specify a pricing scheme in order to define the efficient quantity and welfare. However, the efficient outcome can be achieved if the market intermediary sets a price of $\max \{C_{(K)},V_{[K+1]}\}$ for buyers and $\min\{C_{(K+1)},V_{[K]}\}$ for sellers, where we set $V_{[N+1]} = a$ and $C_{(M+1)} = d$ for convenience (see Loertscher and Marx~\cite{lm13}). Furthermore, the \textquotedblleft matching" of buyers and sellers that appears in \eqref{eq:kdef} and \eqref{eq:welfQ} does not necessarily indicate that these agents trade directly with one another. Indeed, because of the homogeneous goods assumption, it does not matter which buyer trades with which seller. The \textquotedblleft matching" in \eqref{eq:kdef} and \eqref{eq:welfQ} is an algorithm for computing the efficient quantity traded and welfare.}

\pgfplotsset{
    standard/.style={
        xmin=0,xmax=1.3,
        ymin=0,ymax=1.2,
        axis x line=bottom,
        axis y line=middle,
        every axis x label/.style={at={(1,0)},anchor=north},
        every axis y label/.style={at={(0,1)},anchor=east}
    }
}

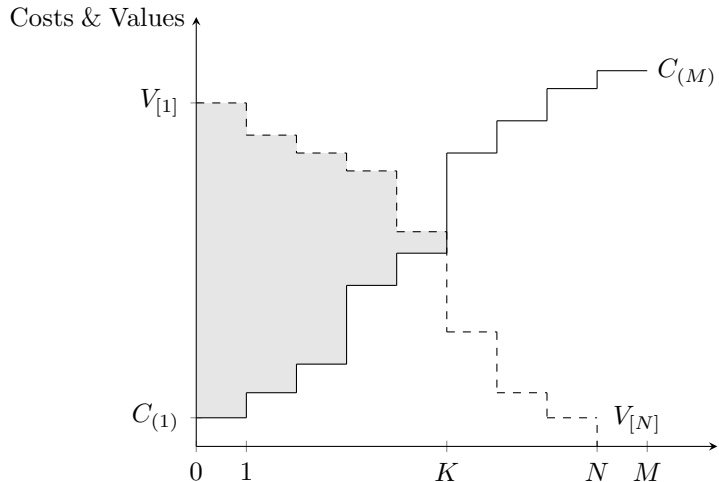
\begin{figure}
	\centering
	\begin{tikzpicture}
	\begin{axis}[
    standard,
    ylabel=\footnotesize{$\mathrm{Costs\; \& \; Values}$},
    xtick={0,0.125,0.625,1,1.125},
    xticklabels={\footnotesize{$0$},\footnotesize{$1$},\footnotesize{$K$},\footnotesize{$N$},\footnotesize{$M$}},
    ytick={0.96,0.08},
    yticklabels={\footnotesize{$V_{[1]}$},\footnotesize{$C_{(1)}$}}
]
		\fill[color=black!10] (axis cs:0.002,0.08) rectangle (axis cs:0.125,0.96);
		\fill[color=black!10] (axis cs:0.125,0.15) rectangle (axis cs:0.25,0.87);
		\fill[color=black!10] (axis cs:0.25,0.23) rectangle (axis cs:0.375,0.82);
		\fill[color=black!10] (axis cs:0.375,0.45) rectangle (axis cs:0.5,0.77);
		\fill[color=black!10] (axis cs:0.5,0.54) rectangle (axis cs:0.625,0.6);
		\addplot[black, no markers] coordinates {(0,0.08) (0.125,0.08)};
		\addplot[black, no markers] coordinates {(0.125,0.08) (0.125,0.15)};
		\addplot[black, no markers] coordinates {(0.125,0.15) (0.25,0.15)};
		\addplot[black, no markers] coordinates {(0.25,0.15) (0.25,0.23)};
		\addplot[black, no markers] coordinates {(0.25,0.23) (0.375,0.23)};
		\addplot[black, no markers] coordinates {(0.375,0.23) (0.375,0.45)};
		\addplot[black, no markers] coordinates {(0.375,0.45) (0.5,0.45)};
		\addplot[black, no markers] coordinates {(0.5,0.45) (0.5,0.54)};
		\addplot[black, no markers] coordinates {(0.5,0.54) (0.625,0.54)};
		\addplot[black, no markers] coordinates {(0.625,0.54) (0.625,0.82)};
		\addplot[black, no markers] coordinates {(0.625,0.82) (0.75,0.82)};
		\addplot[black, no markers] coordinates {(0.75,0.82) (0.75,0.91)};
		\addplot[black, no markers] coordinates {(0.75,0.91) (0.875,0.91)};
		\addplot[black, no markers] coordinates {(0.875,0.91) (0.875,1)};
		\addplot[black, no markers] coordinates {(0.875,1) (1,1)};
        \addplot[black, no markers] coordinates {(1,1) (1,1.05)};
        \addplot[black, no markers] coordinates {(1,1.05) (1.125,1.05)}
		[xshift=15pt,yshift=-10pt]
        node[above] {\footnotesize{$C_{(M)}$}};
		
		\addplot[black, dashed, no markers] coordinates {(0,0.96) (0.125,0.96)};
		\addplot[black, dashed, no markers] coordinates {(0.125,0.96) (0.125,0.87)};
		\addplot[black, dashed, no markers] coordinates {(0.125,0.87) (0.25,0.87)};
		\addplot[black, dashed, no markers] coordinates {(0.25,0.87) (0.25,0.82)};
		\addplot[black, dashed, no markers] coordinates {(0.25,0.82) (0.375,0.82)};
		\addplot[black, dashed, no markers] coordinates {(0.375,0.82) (0.375,0.77)};
		\addplot[black, dashed, no markers] coordinates {(0.375,0.77) (0.5,0.77)};
		\addplot[black, dashed, no markers] coordinates {(0.5,0.77) (0.5,0.6)};
		\addplot[black, dashed, no markers] coordinates {(0.5,0.6) (0.625,0.6)};
		\addplot[black, dashed, no markers] coordinates {(0.625,0.6) (0.625,0.32)};
		\addplot[black, dashed, no markers] coordinates {(0.625,0.32) (0.75,0.32)};
		\addplot[black, dashed, no markers] coordinates {(0.75,0.32) (0.75,0.15)};
		\addplot[black, dashed, no markers] coordinates {(0.75,0.15) (0.875,0.15)};
		\addplot[black, dashed, no markers] coordinates {(0.875,0.15) (0.875,0.08)};
		\addplot[black, dashed, no markers] coordinates {(0.875,0.08) (1,0.08)}
		[xshift=15pt,yshift=-10pt]
        node[above] {\footnotesize{$V_{[N]}$}};
		\addplot[black, dashed, no markers] coordinates {(1,0.06) (1,0)};
\end{axis}
\end{tikzpicture}	
\caption{\small The efficient quantity $K$ is given by the abscissa of the intersection of the   plots of the buyer and seller order statistics, and the respective value of welfare $W$ is equal to the area of the shaded region.}
\label{sd}
\end{figure}

Before stating the main results of this paper, we must introduce some additional model assumptions. First of all, in our large market setup we will consider a family of markets of increasing size with the same distribution functions $F$ and $G$ that satisfy the following standard mechanism design assumption:
\medskip

{\bf(A1)}~{\em The distribution functions $F$ and $G$ are absolutely continuous, with respective densities $f(x)$ and $g(x)$ bounded and bounded away from zero on their respective supports  $[a,b]$ and $[c,d]$, such that $(a,b)\cap (c,d)\neq \varnothing$ $($in other words, $a < d$ and $c < b)$.}
\medskip

\pu{Given the distribution functions $F$ and $G$,} it is natural and convenient to consider a net \pu{(also known as a Moore--Smith sequence)} of markets
\begin{align}
\notag \mathfrak{M}_\alpha
 := \langle \alpha := (M,N),F,G \rangle, \quad \alpha \in \mathscr{A},
\end{align}
indexed by the directed set~$\mathscr{A}$ about which we will make the following assumption:

\medskip
{\bf (A2)} {\em We assume that
\begin{align*}
\mathscr{A} := \{ \alpha  = (M,N) \in \mathbb{N}^2 :
                  \lambda_\alpha := M N^{-1} \in  I\},
\end{align*}
where
$I:=\left[1 - F(d) + \epsilon, 1/(G(a) + \epsilon) \right]$
for a fixed   $\epsilon > 0.$ }

\medskip

Here $\epsilon$ is assumed to be small enough so that $I\neq\varnothing$; note that $F(d) >0$ and $G(a) <1$ by virtue of the assumption $a<d$, see~{\bf (A1)}. Note also that $\inf I >0.$

The set $\mathscr{A} $ is endowed with the natural  preorder: for  $\alpha=(N,M) $ and $\alpha'=(N',M'),$ one has  $\alpha\le \alpha'$ iff $N\le N'$ and $M\le M'.$ We will be interested in the limiting distributions of the efficient quantities  $K_\alpha$ and welfares $W_\alpha$ for the respective markets from the net $\{\mathfrak{M}_\alpha\}_{\alpha \mathscr{A}}$.

The assumption  $\lambda_\alpha\in I$ from {\bf (A2)} excludes trivial cases.  Indeed, with probability tending to one, for large $\alpha$ the number of $V_i$'s exceeding $d$ will be equal to $(1-F(d)+o(1))N$. So if   $\lambda_\alpha < 1 - F(d)-\epsilon$ for a fixed $\epsilon >0$, it would mean that the total number of sellers $M=\lambda_\alpha N$ would be less than the number of buyers with valuations higher than the maximum possible production cost~$d$, meaning that all sellers trade and there is rationing on the demand side of the market. Likewise, the situation when $\lambda_\alpha >  1/G(a)$ corresponds to a market with excess supply (all buyers trade and there is rationing on the supply side of the market).

To state the main results, we need some further notation. We will frequently deal with scaled functions of the form $h(\lambda^{-1}_{\alpha} t)$. For convenience, for any $h:[0,1]\to \mathbb{R}$  we define
\begin{align*}
	  \widehat{h}(s) := h(s \wedge 1),\quad s \geq 0 .
\end{align*}
Note that the function $\widehat{h}(\lambda^{-1}_\alpha t),$ $t \in [0,1],$ is   well-defined   for all $\lambda_\alpha>0$.

Using notation $h^{(-1)}$ for the inverse of function $h$ (to avoid confusion with the reciprocal~$h^{-1}$), introduce the functions
\begin{align}
   \label{eq:ENdef}
   E_\alpha(t) := F^{(-1)} \left( 1-t \right) - \widehat{G^{(-1)}} ( \lambda^{-1}_\alpha t ),\quad t \in [0,1],
\end{align}
and put
\begin{align}
	\label{eq:Hdef}
\mathcal{H}(h) := \sup\{ t \in
	(0,1) : h(t) \geq 0 \}
\end{align}
(which is well-defined for any $h:[0,1]\to \mathbb{R}$ with $h(0) > 0$) and
\begin{align}
	\label{eq:t0def}
 t_{\alpha} := \mathcal{H}(E_{\alpha}) \in (0,\lambda_\alpha \wedge 1),
 \end{align}
where the right-hand relation holds due to the assumption on $I$ from~{\bf (A2)}.

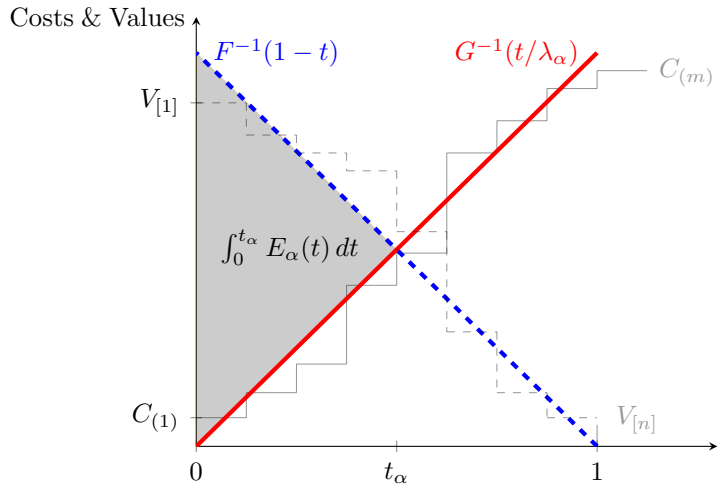
\begin{figure}
	\centering
	\begin{tikzpicture}
	\begin{axis}[
    standard,
    ylabel=\footnotesize{$\mathrm{Costs\; \& \; Values}$},
    xtick={0,0.5,1},
    xticklabels={\footnotesize{$0$},\footnotesize{$t_\alpha$},\footnotesize{$1$}},
    ytick={0.96,0.08},
    yticklabels={\footnotesize{$V_{[1]}$},\footnotesize{$C_{(1)}$}}
]
		\addplot[black, no markers, opacity=0.4] coordinates {(0,0.08) (0.125,0.08)};
		\addplot[black, no markers, opacity=0.4] coordinates {(0.125,0.08) (0.125,0.15)};
		\addplot[black, no markers, opacity=0.4] coordinates {(0.125,0.15) (0.25,0.15)};
		\addplot[black, no markers, opacity=0.4] coordinates {(0.25,0.15) (0.25,0.23)};
		\addplot[black, no markers, opacity=0.4] coordinates {(0.25,0.23) (0.375,0.23)};
		\addplot[black, no markers, opacity=0.4] coordinates {(0.375,0.23) (0.375,0.45)};
		\addplot[black, no markers, opacity=0.4] coordinates {(0.375,0.45) (0.5,0.45)};
		\addplot[black, no markers, opacity=0.4] coordinates {(0.5,0.45) (0.5,0.54)};
		\addplot[black, no markers, opacity=0.4] coordinates {(0.5,0.54) (0.625,0.54)};
		\addplot[black, no markers, opacity=0.4] coordinates {(0.625,0.54) (0.625,0.82)};
		\addplot[black, no markers, opacity=0.4] coordinates {(0.625,0.82) (0.75,0.82)};
		\addplot[black, no markers, opacity=0.4] coordinates {(0.75,0.82) (0.75,0.91)};
		\addplot[black, no markers, opacity=0.4] coordinates {(0.75,0.91) (0.875,0.91)};
		\addplot[black, no markers, opacity=0.4] coordinates {(0.875,0.91) (0.875,1)};
		\addplot[black, no markers, opacity=0.4] coordinates {(0.875,1) (1,1)};
        \addplot[black, no markers, opacity=0.4] coordinates {(1,1) (1,1.05)};
        \addplot[black, no markers, opacity=0.4] coordinates {(1,1.05) (1.125,1.05)}
		[xshift=15pt,yshift=-10pt]
        node[above] {\footnotesize{$C_{(m)}$}};
		
		\addplot[patch, black, opacity=0.2]
    	coordinates {(0,0) (0,1.1) (0.5,0.55)};
		\addplot[ultra thick, blue, dashed, no markers] coordinates {(1,0) (0,1.1)}
		[xshift=30pt,yshift=-10pt]
        node[above] {\footnotesize{$F^{-1}(1-t)$}};
		\addplot[ultra thick, red, no markers] coordinates {(0,0) (1,1.1)}
		[xshift=-30pt,yshift=-10pt]
        node[above] {\footnotesize{$G^{-1}(t/\lambda_\alpha)$}};
		\addplot[black, no markers] coordinates {(0.5,0.55) (0.5,0.55)}
		[xshift=-40pt,yshift=-10pt]
        node[above] {\footnotesize{$\int_{0}^{t_\alpha} E_\alpha(t) \, dt $}};

		\addplot[black, dashed, no markers, opacity=0.4] coordinates {(0,0.96) (0.125,0.96)};
		\addplot[black, dashed, no markers, opacity=0.4] coordinates {(0.125,0.96) (0.125,0.87)};
		\addplot[black, dashed, no markers, opacity=0.4] coordinates {(0.125,0.87) (0.25,0.87)};
		\addplot[black, dashed, no markers, opacity=0.4] coordinates {(0.25,0.87) (0.25,0.82)};
		\addplot[black, dashed, no markers, opacity=0.4] coordinates {(0.25,0.82) (0.375,0.82)};
		\addplot[black, dashed, no markers, opacity=0.4] coordinates {(0.375,0.82) (0.375,0.77)};
		\addplot[black, dashed, no markers, opacity=0.4] coordinates {(0.375,0.77) (0.5,0.77)};
		\addplot[black, dashed, no markers, opacity=0.4] coordinates {(0.5,0.77) (0.5,0.6)};
		\addplot[black, dashed, no markers, opacity=0.4] coordinates {(0.5,0.6) (0.625,0.6)};
		\addplot[black, dashed, no markers, opacity=0.4] coordinates {(0.625,0.6) (0.625,0.32)};
		\addplot[black, dashed, no markers, opacity=0.4] coordinates {(0.625,0.32) (0.75,0.32)};
		\addplot[black, dashed, no markers, opacity=0.4] coordinates {(0.75,0.32) (0.75,0.15)};
		\addplot[black, dashed, no markers, opacity=0.4] coordinates {(0.75,0.15) (0.875,0.15)};
		\addplot[black, dashed, no markers, opacity=0.4] coordinates {(0.875,0.15) (0.875,0.08)};
		\addplot[black, dashed, no markers, opacity=0.4] coordinates {(0.875,0.08) (1,0.08)}
		[xshift=15pt,yshift=-10pt]
        node[above] {\footnotesize{$V_{[n]}$}};
		\addplot[black, no markers, opacity=0.4] coordinates {(1,0.06) (1,0)};
\end{axis}
\end{tikzpicture} 
\caption{\small The scaled empirical quantile functions, which are associated with demand and supply, may be approximated by the associated quantile functions. The errors associated with these approximations are of order $1/\sqrt{N}$ and are given by appropriately scaled Brownian bridges.}
\label{fig:approx}	
\end{figure}

\begin{figure}
    \centering
    \begin{subfigure}[b]{0.32\textwidth}
        \includegraphics[width=\textwidth]{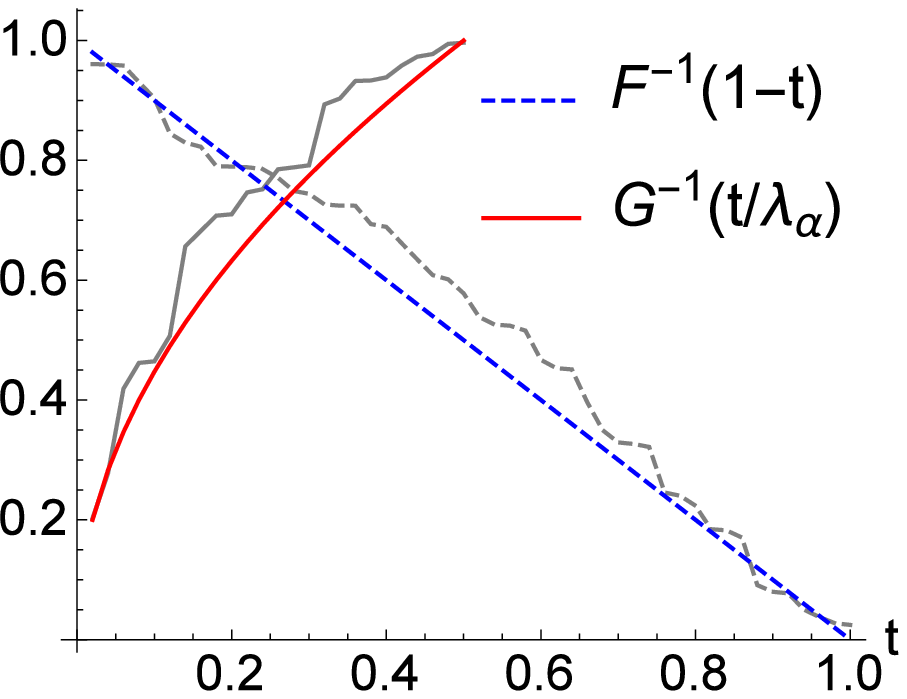}
         \caption{$\alpha = (50,25)$}
    \end{subfigure}
    \begin{subfigure}[b]{0.32\textwidth}
        \includegraphics[width=\textwidth]{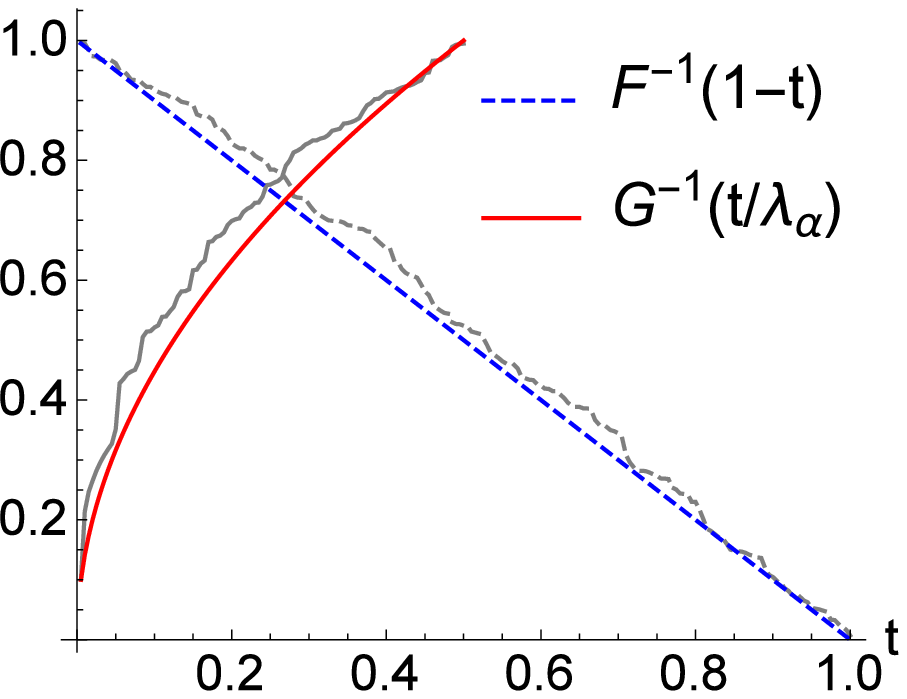}
        \caption{$\alpha = (200,100)$}
    \end{subfigure}
    \begin{subfigure}[b]{0.32\textwidth}
        \includegraphics[width=\textwidth]{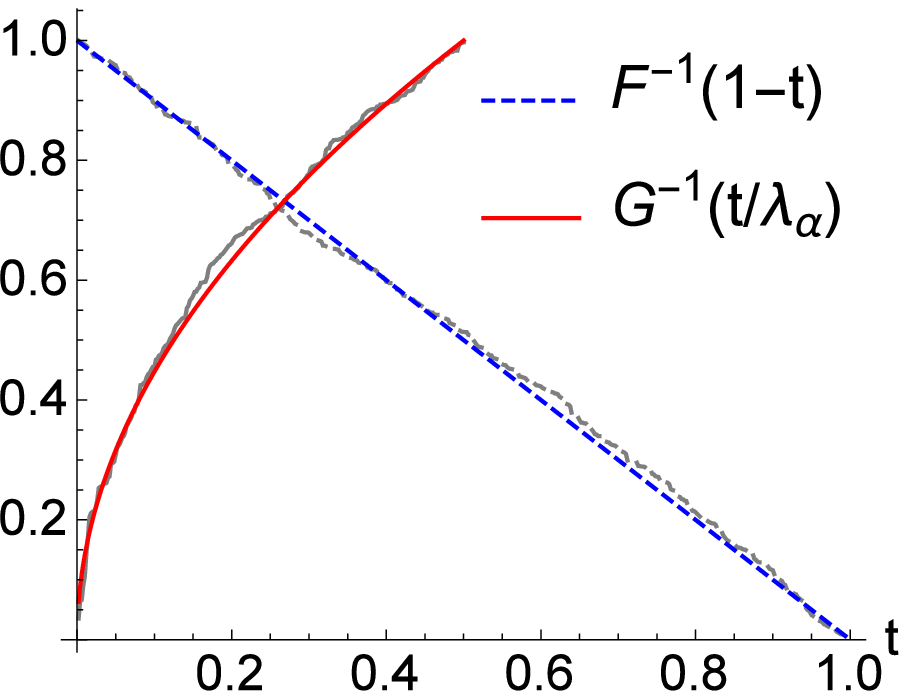}
        \caption{$\alpha = (500,250)$}
    \end{subfigure}
    \caption{\pu{Three simulations illustrating the approximation of the scaled emirical quantile functions by the theoretical quantile functions when $F(x) =x$ and  $G(x)=x^2$ for $x\in [0,1]$.}}
    \label{fig:sim}
\end{figure}

Due to the a.s.\ convergence of empirical quantile functions to the theoretical ones, for large markets the function $E_\alpha (t)$ approximates the difference between the step-functions whose plots are depicted, respectively, by the dashed and solid lines in Figure~\ref{fig:approx} (cf.~\eqref{eq:Kalpha} below). \pu{Simulations illustrating this approximation are shown in Figure \ref{fig:sim}.} So $t_\alpha N$ will be the ``first order approximation'' to~$K_\alpha$, while the integral of $E_\alpha$ over $(0,t_\alpha)$ will specify a deterministic approximation to~$W_\alpha$. The ``second order approximation" to both $K_\alpha$ and $W_\alpha$ will be obtained in this paper using the second order approximation to empirical quantile functions provided by the sum of the theoretical quantile function and a Brownian bridge process.

Now observe that, provided that $f,g$ are continuous inside their respective supports, the function $E_\alpha(t)$ is continuously differentiable for $t \in (0,1) \setminus \{  \lambda_\alpha\}$, and
\begin{align}
	\label{eq:ENderiv} E_{\alpha}'(t) = - \frac{1}{f(F^{(-1)}(1-t))} - \frac{1}{\lambda_\alpha g(G^{(-1)}(\lambda_\alpha^{-1}t))}, \quad t \in (0,\lambda_\alpha \wedge 1).
\end{align}

We will  need one more   technical assumption  on the distribution functions $F$ and~$G$:
\medskip

{\bf(A3)}~{\em The densities $f$ and $g$ are differentiable on $(a,b)$ and $(c,d),$ respectively. Moreover, the functions
\[
\frac{d}{dt} \frac{1}{f(F^{(-1)}( t))}
 = -\frac{f'(F^{(-1)}(t))}{f^3(F^{(-1)}(t))}
 \quad \mbox{and}\quad
\frac{d}{dt} \frac{1}{g (G^{(-1)}( t))}
 = - \frac{g'(G^{(-1)}(t))}{g^3(G^{(-1)}(t))}
\]
are bounded on $ (0,1)$.}
\medskip

Finally, let
 \begin{align}
 \label{eq:sigmadef}
	  \sigma_\alpha^2
    : = \frac{1}{(E_\alpha'(t_\alpha))^2} \left[ \frac{t_{\alpha}(1-t_{\alpha})}{f^2(F^{(-1)}(1-t_{\alpha}))}
    +
    \frac{t_{\alpha}(1-\lambda^{-1}_\alpha t_{\alpha})}{\lambda^{2}_\alpha g^2(G^{(-1)}(\lambda^{-1}_\alpha t_{\alpha}))}
    \right],
\end{align}
and
 \begin{align}
 \label{eq:varsigmadef}
 \varsigma_\alpha^2: = 2 \int_0^{t_\alpha} S_\alpha (t) dt ,
 \end{align}
 where the function $S_\alpha (t) ,$ $t\in (0,t_\alpha),$ is given by
\begin{align*}
 \frac{1-t}{f(F^{(-1)}(1-t))}\int_{F^{(-1)}(1-t)}^{b}  (1-F(x))  dx
 + \frac{ 1 - \lambda^{-1}_\alpha t}{ g(G^{(-1)}(\lambda^{-1}_\alpha t))}
 \int_c^{G^{(-1)}(\lambda^{-1}_\alpha t )} G(x)  dx     .
\end{align*}

Our main result is the following strong approximation theorem. It implies  that, under the above assumptions, the efficient quantity and welfare are asymptotically normal for large markets. \pu{More precisely, the univariate distributions of $K_\alpha$ and $W_\alpha$ can be approximated by normal distributions with respective means and variances $(N t_\alpha, N \sigma_\alpha^2)$ and $\bigl(N \int_{0}^{t_\alpha} E_\alpha (t)\, dt, N \varsigma_\alpha^2\bigr)$. Moreover, the joint distribution of $(K_\alpha,W_\alpha)$ is asymptotically normal as well.} We also give upper bounds for the convergence rates.

\begin{theorem}
\label{thm:1}
Under assumptions {\bf (A1)}--{\bf (A3)}, for the net of markets $\{\mathfrak{M}_\alpha\}_{\alpha\in\mathscr{A}}$   there exist a net $\{\mathfrak{R}_\alpha\}_{\alpha\in\mathscr{A}}$ of realizations of these markets on a common probability space together with a net of bivariate normal random vectors $\{(Z^{(1)}_\alpha, Z^{(2)}_\alpha)\}_{\alpha\in\mathscr{A}}$ with zero means and
\[
{\rm Var}\,(Z^{(1)}_\alpha)=\sigma_\alpha^2, \quad
 {\rm Var}\,(Z^{(2)}_\alpha)=\varsigma_\alpha^2, \quad
 {\rm Cov}\,(Z^{(1)}_\alpha,Z^{(2)}_\alpha)= \varkappa_\alpha
  := -S_\alpha (t_\alpha)/E'(t_\alpha) ,
\]
such that \pu{for the efficient quantity $K_\alpha$ one has}
\begin{align}
 \label{eq:N1/4}
	  \limsup_\alpha  \frac{N^{1/4}}{(\ln N)^{1/2}}\left| \frac{K_\alpha -  N t_\alpha  }{N^{1/2}  }  -Z^{(1)}_\alpha \right| <\infty \quad \mbox{a.s.}
\end{align}
and \pu{for welfare $W_\alpha$ one has}
\begin{align*}
 \limsup_\alpha  \frac{N^{1/2}}{ \ln N }\left| \frac{1}{N^{1/2} }
 \left(W_\alpha - N \int_{0}^{t_\alpha} E_\alpha (t)\, dt \right)  -Z^{(2)}_\alpha \right| <\infty \quad \mbox{a.s.}
\end{align*}	
\end{theorem}

\pu{The proof uses the common probability space method and is deferred to Section~\ref{sec:2}. The idea is to consider, for a given $\alpha = (N,M),$ the empirical quantile functions associated with the samples of buyer valuations $V_1, \dots, V_N$ and seller costs $C_1,\dots,C_M$.  The efficient quantity and welfare, defined in \eqref{eq:kdef} and \eqref{eq:welfQ} respectively, can then be expressed as linear functionals of the respective  empirical quantile processes. By the results of Cs\"{o}rg\H{o} and R\'{e}v\'{e}sz~\cite{cr75}, it is possible to construct a suitable probability space carrying a sequence of empirical quantile processes that converge almost surely to an appropriately weighted Brownian bridge process. We use this result to approximate $K_\alpha$ and $W_\alpha$ by linear functionals of the respective limiting Gaussian process. To complete the proof, we derive and apply an appropriate generalisation of the delta method and compute bounds for the approximation rate.}

\begin{figure}[ht]
	\centering
	{\includegraphics[scale=1]{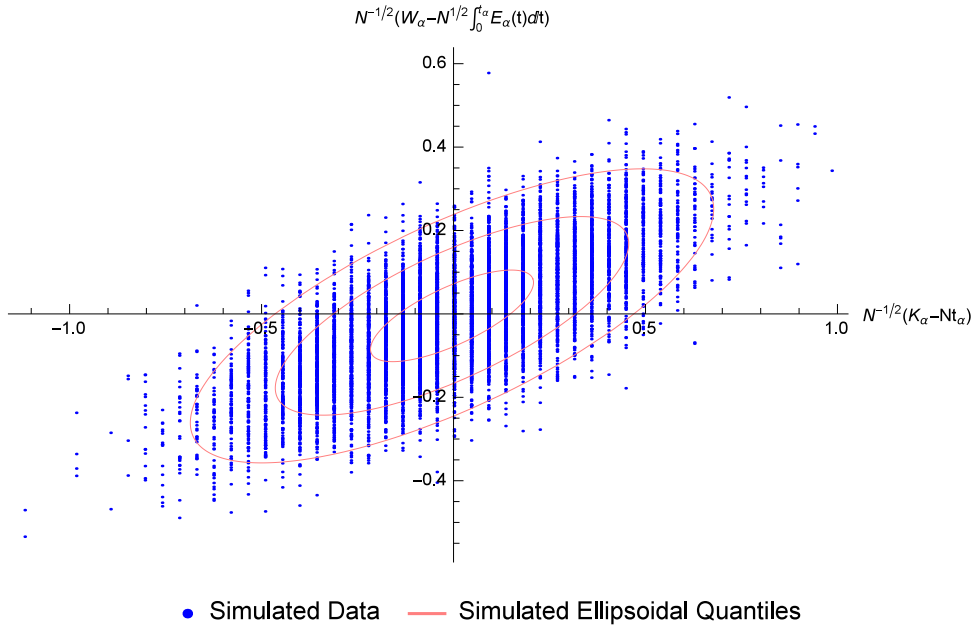}}	
	 	\caption{\small The scatterplot of simulated $10^4$ independent scaled copies of  $(K_\alpha, W_\alpha)$ and   its ellipsoidal quantiles (see the paragraph  following Theorem~\ref{thm:1} for more detail).}
 \label{scatter}
\end{figure}

\begin{figure}[ht]
	\centering{\includegraphics[scale=1]{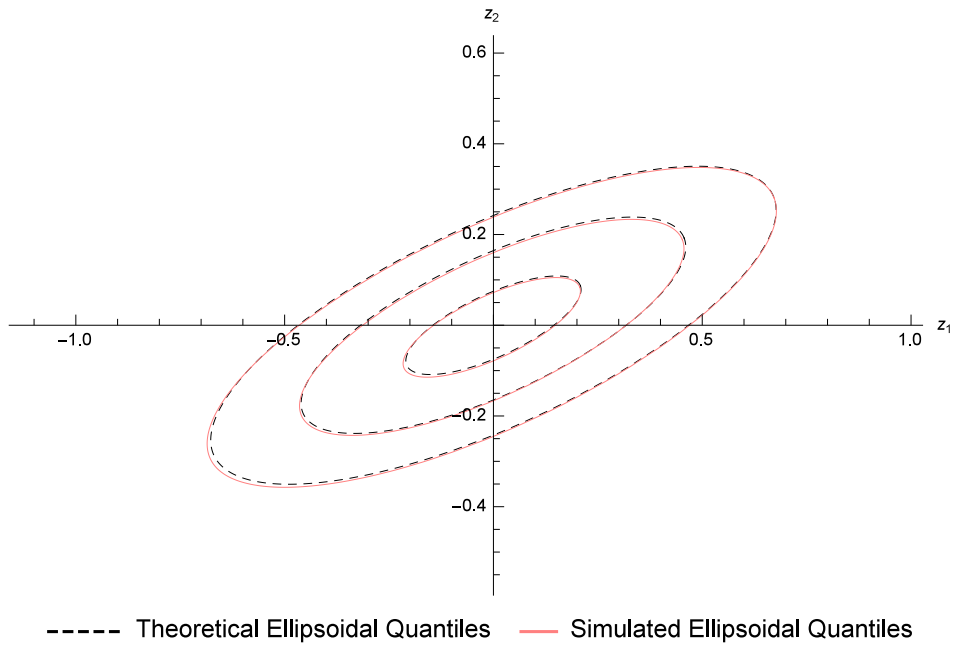}}
\caption{\small The empirical ellipsoidal quantiles from Figure~\ref{scatter} together with the ellipsoidal quantiles of the approximating bivariate normal distribution.}
 \label{ellipse}
\end{figure}

As an illustration to the assertion of the theorem, Figure~\ref{scatter} shows the scatterplot of $10^4$ independent realizations of the centered and scaled (as  per the statement of Theorem~\ref{thm:1}) random  vector  $(K_\alpha, W_\alpha)$ in the case when $\alpha=(500,250),$ $F(x) =x$ and  $G(x)=x^2$ for $x\in [0,1]$, together with the 0.25, 0.5 and 0.95 ellipsoidal quantiles for the sample. Figure~\ref{ellipse} displays these empirical ellipsoidal quantiles together with the respective (theoretical) ellipsoidal quantiles for the approximating distribution of~$(Z^{(1)}_\alpha, Z^{(2)}_\alpha)$. Both plots were generated using {\sc Mathematica~10}.

\pu{Recall that, for a bivariate distribution $P$ with finite second order moments, the ellipsoidal quantiles are defined as follows. Let $\mu$ and $\Sigma$ be the mean vector and covariance matrix of $P$, respectively. The ellipsoidal quantile of level $u \in (0,1)$ is defined as the boundary $\partial A$ of the smallest set of the form}
\begin{align*}
	\pu{A := \{y \in \mathbb{R}^2:(y - \mu) \Sigma (y-\mu)^T \leq \mathrm{const}\}}
\end{align*}
\pu{such that $P(A) \geq u$. Empirical ellipsoidal quantiles for a bivariate sample are defined as the ellipsoidal quantiles of the respective empirical distribution.}

\begin{remark}
{\rm
Note that $N\asymp M$ ($N$ and $M$ are asymptotically equivalent) under assumption {\bf (A2)}, so one could state the assertions of the theorem in a similar way using $M$ rather than $N$ as well. Note also that $\varkappa_\alpha>0$  (as one would expect, of course) since $S_\alpha (t_\alpha)>0$ and $E_\alpha' (t_\alpha)<0.$
}
\end{remark}

\begin{remark}
{\rm
The true rate of convergence of the distribution of $K_\alpha$ to the normal distribution is most likely $N^{-1/2}$, as indicated by the bounds~\eqref{eq:BerryEsseen} below established in the special case when~$F=G$. So the crudeness of the bound in~\eqref{eq:N1/4} seems to be due to the method of proof employed.
}
\end{remark}

In the interesting special case $F=G,$ the distribution of~$K_\alpha$ is independent of~$F$. Moreover, it is known to be hypergeometric (see the proof of part~(ii) below for a more precise statement) and the parameters of the approximating normal law for~$K_\alpha$ admit   simple  explicit representations as  functions of~$\alpha$. In addition, in that case one can establish a better convergence rate to the normal distribution than that claimed in Theorem~\ref{thm:1}.

The distribution of $W_\alpha$ (and hence that of $(K_\alpha, W_\alpha)$) and the approximating normal law do depend on~$F$. However,  one can also obtain simple closed formulae for the parameters of the approximating normal laws when the common distribution $F=G$ is uniform on $(a,b)$. We will only deal here with the case $a=0,$ $b=1$, the results in the general case following in a straightforward way, and state our results in the form of the following theorem.

\begin{theorem}
 \label{thm:2}
If $F = G$ then:
\medskip

{\rm (i)}~one has
\begin{align*}
	t_\alpha = \frac{\lambda_\alpha}{1 + \lambda_\alpha}, \qquad   \sigma_\alpha^2 = \frac{\lambda_\alpha^2}{(1+\lambda_\alpha)^3};
\end{align*}

{\rm (ii)}~there exist constants $C_1, C_2\in (0,\infty)$ such that
\begin{align}
 \label{eq:BerryEsseen}
\frac{C_1}{N^{1/2}}
 \le
 \sup_{x}\biggl| {\bf P}  \left(\frac{K_\alpha - N t_\alpha}{\sigma_\alpha N^{1/2}}\le x \right) - \Phi (x) \biggr|
 \le
\frac{C_2}{N^{1/2}},
\end{align}
where $\Phi$ is  the standard normal distribution function.

{\rm (iii)} If, moreover,    $F(t) = G(t) = t$ for $t\in (0,1),$ then
\begin{align}
\label{eq:varsigma}
	  \int_{0}^{t_\alpha} E_\alpha(t) \, dt = \frac{\lambda_\alpha}{2(1+\lambda_\alpha)},
 \quad
 \varsigma_\alpha^2 = \frac{\lambda_\alpha(1+3\lambda_\alpha + \lambda_\alpha^2)}{12(1+\lambda_\alpha)^3},
 \quad
\varkappa_\alpha  =  \frac{\lambda_\alpha^2}{2(1+\lambda_\alpha)^3}.
\end{align}	
\end{theorem}

\begin{remark}
{\rm
Observe that, in the case from part~(iii), the correlation coefficient between the components of the approximating normal distribution can be easily found to be equal to $\sqrt{3  /(\lambda_\alpha^{-1} +3 +\lambda_\alpha )}.$ This quantity  attains its maximum value $\sqrt{3/5}$ at $\lambda_\alpha=1 $ (that is, when $M=N$) and vanishes as $\lambda_\alpha  \vee \lambda_\alpha^{-1} \to \infty. $}
\end{remark}

\section{Applications and extensions}

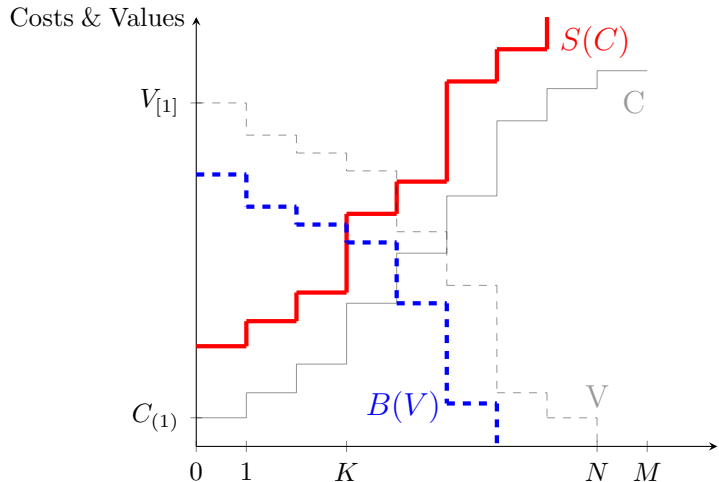
\begin{figure}
	\centering
	\begin{tikzpicture}
	\begin{axis}[
    standard,
    ylabel=\footnotesize{$\mathrm{Costs\; \& \; Values}$},
    xtick={0,0.125,0.375,1,1.125},
    xticklabels={\footnotesize{$0$},\footnotesize{$1$},\footnotesize{$K$},\footnotesize{$N$},\footnotesize{$M$}},
    ytick={0.96,0.08},
    yticklabels={\footnotesize{$V_{[1]}$},\footnotesize{$C_{(1)}$}}
]
		\addplot[black, no markers, opacity=0.4] coordinates {(0,0.08) (0.125,0.08)};
		\addplot[black, no markers, opacity=0.4] coordinates {(0.125,0.08) (0.125,0.15)};
		\addplot[black, no markers, opacity=0.4] coordinates {(0.125,0.15) (0.25,0.15)};
		\addplot[black, no markers, opacity=0.4] coordinates {(0.25,0.15) (0.25,0.23)};
		\addplot[black, no markers, opacity=0.4] coordinates {(0.25,0.23) (0.375,0.23)};
		\addplot[black, no markers, opacity=0.4] coordinates {(0.375,0.23) (0.375,0.4)};
		\addplot[black, no markers, opacity=0.4] coordinates {(0.375,0.4) (0.5,0.4)};
		\addplot[black, no markers, opacity=0.4] coordinates {(0.5,0.4) (0.5,0.54)};
		\addplot[black, no markers, opacity=0.4] coordinates {(0.5,0.54) (0.625,0.54)};
		\addplot[black, no markers, opacity=0.4] coordinates {(0.625,0.54) (0.625,0.7)};
		\addplot[black, no markers, opacity=0.4] coordinates {(0.625,0.7) (0.75,0.7)};
		\addplot[black, no markers, opacity=0.4] coordinates {(0.75,0.7) (0.75,0.91)};
		\addplot[black, no markers, opacity=0.4] coordinates {(0.75,0.91) (0.875,0.91)};
		\addplot[black, no markers, opacity=0.4] coordinates {(0.875,0.91) (0.875,1)};
		\addplot[black, no markers, opacity=0.4] coordinates {(0.875,1) (1,1)};
		\addplot[black, no markers, opacity=0.4] coordinates {(1,1) (1,1.05)};
		\addplot[black, no markers, opacity=0.4] coordinates {(1,1.05) (1.125,1.05)}
		[xshift=-5pt,yshift=-20pt]
        node[above] {C};
        
		\addplot[ultra thick, red, no markers] coordinates {(0,0.28) (0.125,0.28)};
		\addplot[ultra thick, red, no markers] coordinates {(0.125,0.28) (0.125,0.35)};
		\addplot[ultra thick, red, no markers] coordinates {(0.125,0.35) (0.25,0.35)};
		\addplot[ultra thick, red, no markers] coordinates {(0.25,0.35) (0.25,0.43)};
		\addplot[ultra thick, red, no markers] coordinates {(0.25,0.43) (0.375,0.43)};
		\addplot[ultra thick, red, no markers] coordinates {(0.375,0.43) (0.375,0.65)};
		\addplot[ultra thick, red, no markers] coordinates {(0.375,0.65) (0.5,0.65)};
		\addplot[ultra thick, red, no markers] coordinates {(0.5,0.65) (0.5,0.74)};
		\addplot[ultra thick, red, no markers] coordinates {(0.5,0.74) (0.625,0.74)};
		\addplot[ultra thick, red, no markers] coordinates {(0.625,0.74) (0.625,1.02)};
		\addplot[ultra thick, red, no markers] coordinates {(0.625,1.02) (0.75,1.02)};
		\addplot[ultra thick, red, no markers] coordinates {(0.75,1.02) (0.75,1.11)};
		\addplot[ultra thick, red, no markers] coordinates {(0.75,1.11) (0.875,1.11)};
		\addplot[ultra thick, red, no markers] coordinates {(0.875,1.11) (0.875,1.2)}			
		[xshift=18pt,yshift=-20pt]
        node[above] {$S(C)$};

		\addplot[black, dashed, no markers, opacity=0.4] coordinates {(0,0.96) (0.125,0.96)};
		\addplot[black, dashed, no markers, opacity=0.4] coordinates {(0.125,0.96) (0.125,0.87)};
		\addplot[black, dashed, no markers, opacity=0.4] coordinates {(0.125,0.87) (0.25,0.87)};
		\addplot[black, dashed, no markers, opacity=0.4] coordinates {(0.25,0.87) (0.25,0.82)};
		\addplot[black, dashed, no markers, opacity=0.4] coordinates {(0.25,0.82) (0.375,0.82)};
		\addplot[black, dashed, no markers, opacity=0.4] coordinates {(0.375,0.82) (0.375,0.77)};
		\addplot[black, dashed, no markers, opacity=0.4] coordinates {(0.375,0.77) (0.5,0.77)};
		\addplot[black, dashed, no markers, opacity=0.4] coordinates {(0.5,0.77) (0.5,0.6)};
		\addplot[black, dashed, no markers, opacity=0.4] coordinates {(0.5,0.6) (0.625,0.6)};
		\addplot[black, dashed, no markers, opacity=0.4] coordinates {(0.625,0.6) (0.625,0.45)};
		\addplot[black, dashed, no markers, opacity=0.4] coordinates {(0.625,0.45) (0.75,0.45)};
		\addplot[black, dashed, no markers, opacity=0.4] coordinates {(0.75,0.45) (0.75,0.15)};
		\addplot[black, dashed, no markers, opacity=0.4] coordinates {(0.75,0.15) (0.875,0.15)};
		\addplot[black, dashed, no markers, opacity=0.4] coordinates {(0.875,0.15) (0.875,0.08)};
		\addplot[black, dashed, no markers, opacity=0.4] coordinates {(0.875,0.08) (1,0.08)}
		[xshift=0pt,yshift=0pt]
        node[black,above] {V};
		\addplot[black, dashed, no markers, opacity=0.4] coordinates {(1,0.08) (1,0)};
		
		\addplot[ultra thick, blue, dashed, no markers] coordinates {(0,0.76) (0.125,0.76)};
		\addplot[ultra thick, blue, dashed, no markers] coordinates {(0.125,0.76) (0.125,0.67)};
		\addplot[ultra thick, blue, dashed, no markers] coordinates {(0.125,0.67) (0.25,0.67)};
		\addplot[ultra thick, blue, dashed, no markers] coordinates {(0.25,0.67) (0.25,0.62)};
		\addplot[ultra thick, blue, dashed, no markers] coordinates {(0.25,0.62) (0.375,0.62)};
		\addplot[ultra thick, blue, dashed, no markers] coordinates {(0.375,0.62) (0.375,0.57)};
		\addplot[ultra thick, blue, dashed, no markers] coordinates {(0.375,0.57) (0.5,0.57)};
		\addplot[ultra thick, blue, dashed, no markers] coordinates {(0.5,0.57) (0.5,0.4)};
		\addplot[ultra thick, blue, dashed, no markers] coordinates {(0.5,0.4) (0.625,0.4)};
		\addplot[ultra thick, blue, dashed, no markers] coordinates {(0.625,0.4) (0.625,0.12)};
		\addplot[ultra thick, blue, dashed, no markers] coordinates {(0.625,0.12) (0.75,0.12)};
		\addplot[ultra thick, blue, dashed, no markers] coordinates {(0.75,0.12) (0.75,0)}		[xshift=-35pt,yshift=4pt] node[above] {$B(V)$};	
\end{axis}
\end{tikzpicture}	
\caption{\small Our analysis applies to any mechanism that induces the efficient outcome with respect to non-decreasing transformations of the empirical quantile functions associated with buyer valuations and seller costs.}
\label{transform}
\end{figure}

\pu{For ease of exposition, we restricted attention to the efficient mechanism in Section 1. However, our analysis generalises to a much richer class of mechanisms. Indeed, let $B$ and $S$ be non-decreasing real functions and consider the transformed sets of buyer valuations and seller costs given by $\{B(V_1),\dots,B(V_N)\}$ and $\{S(C_1),\dots,S(C_M)\}$, respectively. Then our analysis immediately applies to the mechanism that induces the efficient outcome with respect to these transformed valuations and costs, provided the assumptions {\bf (A1)}--{\bf (A3)} hold for the composite functions $F \circ B^{(-1)}$ and $G \circ S^{(-1)}$. This is illustrated in Figure~\ref{transform}. Loertscher, Muir and Taylor~\cite{lmt16b} show that this setup subsumes a large class of mechanisms, where the functions $B$ and $S$ depend on the objective of the intermediary and the equilibrium bidding strategies of buyers and sellers. For example, our analysis applies to profit-maximisation and constrained efficient mechanisms. Specifically, standard (but fairly involved) mechanism design arguments show that if we set 
\[
B(x) = x - \frac{1-F(x)}{f(x)}, \quad x \in [a,b], \quad \mathrm{and} \quad S(y) = y + \frac{G(y)}{g(y)}, \quad y \in [c,d],
\]
then we obtain the profit-maximization mechanism.}

\pu{For some mechanisms, such as $k$-double auctions (refer to, for example, Chatterjee and Samuelson~\cite{cs83} and Satterthwaite and Williams~\cite{sw89}), the strategic behaviour of agents depends on $N$ and $M$. More precisely, the strategic behaviour of a given agent refers to the magnitude of the difference between the private information of that agent and the information that agent reveals to the intermediary. Our analysis will only apply to such mechanisms if the strategic behaviour of all buyers and sellers vanishes uniformly on the support of the respective distribution functions $F$ and $G$ at a rate faster than $N^{-1} \ln{N}$ and $M^{-1}\ln{M}$, respectively. For example, by the results of Rustichini, Satterthwaite and Williams~\cite{rsw94}, the $k$-double auction satisfies this condition.}

\section{Proofs of the main results}
\label{sec:2}

The proofs will use the common probability space method and transformed uniform empirical quantile functions, so we will begin with introducing   basic notation and recalling some key facts.

For a sample $X_1,\dots,X_n$   with order statistics       $X_{(1)} \leq \dots \leq X_{(n)}$, the respective empirical quantile function $\mathbb{X}_n (t),$ $t\in [0,1],$ is defined by
\begin{align*}
  \mathbb{X}_n(0) := X_{(1)},\quad  \mathbb{X}_n(t) := X_{(i)}, \quad (i-1)n^{-1} < t \leq in^{-1}, \quad i = 1,\dots,n.
\end{align*}
For the empirical quantile function $\mathbb{U}_n(t)$ constructed from a sample of $\mathrm{U}(0,1)$-distributed independent and identically distributed $U_1,\dots,U_n$,
\[
	\mathbb{R}_n(t) := \sqrt{n} \left( \mathbb{U}_n(t) - t\right), \quad t \in [0,1],
\]
denotes the respective uniform quantile process. Recall that both  $\mathbb{U}_n$  and $\mathbb{R}_n$  are random elements of the  Skorokhod space $\mathscr{D}[0,1]$   of c\`{a}dl\`{a}g functions on~$[0,1]$.

By the well-known Donsker's theorem (see, for example, Section 14 of Billingsley~\cite{billingsley99}), as $n \rightarrow \infty$, the distribution of $\mathbb{R}_n$ in $\mathscr{D}[0,1]$ converges weakly to that of the Brownian bridge process~$\mathbb{B}$. A sharp bound on the rate of convergence is given in Cs\"{o}rg\H{o} and R\'{e}v\'{e}sz~\cite{cr75}. A corollary of that result states that there exists a probability space carrying a sequence of processes $\mathbb{R}^*_n \eqd \mathbb{R}_n$, $n = 1,2,\dots,$ and a Brownian bridge process $\mathbb{B}$ such that, as $n \rightarrow \infty$,
\begin{align*}
	  \vnorm{\mathbb{R}^*_n - \mathbb{B}}_{\infty} = O\left(n^{-1/2} \ln{n} \right)\quad \mathrm{a.s.,}
\end{align*}
where $\vnorm{h}_{\infty} := \sup_{t \in [0,1]} |h(t)|,$ $ h \in \mathscr{D}[0,1].
$

Therefore, as $n \rightarrow \infty$,
\begin{align} 	\label{eq:CRresult}
\mathbb{U}_n(t) = t + n^{-1/2} \mathbb{R}_n(t) \eqd t + n^{-1/2} \mathbb{R}^*_n(t) = t + n^{-1/2} \mathbb{B}(t) + \theta_n(t) n^{-1} \ln{n}
\end{align}
for $t \in [0,1]$, where $\vnorm{\theta_n}_{\infty} = O(1) \; \mathrm{a.s.}$

\medskip
\noindent
{\em Proof of Theorem~\ref{thm:1}.} It is easy to see from~\eqref{eq:kdef},  \eqref{eq:Hdef} and \eqref{eq:t0def}  that
\begin{align}
	\label{eq:Kalpha1}
\delta_\alpha:=  K_\alpha N^{-1} -t_\alpha = \mathcal{H}\left(\mathbb{E}_\alpha \right)-\mathcal{H} (E_\alpha),
\end{align}
where
\begin{align}
	\label{eq:Kalpha}
 \mathbb{E}_\alpha(t) := \mathbb{V}_N \left( 1-t \right) - \widehat{\mathbb{C}_M} ( \lambda^{-1}_\alpha t ), \quad t \in [0,1],
\end{align}
and  $\mathbb{V}_N$ and $\mathbb{C}_M$ denote the empirical quantile functions for the samples of buyer valuations $V_1,\dots,V_N$ and seller costs $C_1,\dots,C_M,$ respectively. Now we will analyze the behavior of $\delta_\alpha$ for ``large''~$\alpha$ when market realizations $\mathcal{R}_\alpha$ are defined using the common probability space construction~\eqref{eq:CRresult}.

To that end, we will assume that
\begin{align}
	\label{eq:vsample}
 V_i = V_{i,N}  := F^{(-1)}\left(U^V_{i,N} \right) \quad  \mathrm{and} \quad  C_j = C_{j,M}  := G^{(-1)}\left(U^C_{j,M} \right),
\end{align}
where $\left\{U^V_{1,N},\dots,U^V_{N,N} \right\}_{N \geq 1}$ and $\left\{U^C_{1,M},\dots,U^C_{M,M} \right\}_{M  \geq 1}$ are independent triangular arrays of row-wise independent $\mathrm{U}(0,1)$-distributed random variables, which are defined on a common probability space with Brownian bridges~$\mathbb{B}^V$ and~$\mathbb{B}^C$ that are independent of each other so that, as $N,M \rightarrow \infty$, for the respective empirical quantile functions one has
\begin{align}
	\begin{split}
	\label{eq:array}
\mathbb{U}^{V}_N(t) & = t + N^{-1/2} \mathbb{B}^V(t) + \theta^V_N(t) N^{-1} \ln{N}, \quad t\in[0,1], \\
	\mathbb{U}^C_M(t) & = t + M^{-1/2} \mathbb{B}^C(t) + \theta^C_M(t) M^{-1} \ln{M}, \quad t \in [0,1],
	\end{split}
\end{align}
where $\vnorm{\theta_N^V}_{\infty} = O(1)$  and $\vnorm{\theta_M^C}_{\infty} = O(1) \; \mathrm{a.s.}$

In this construction, the empirical quantile functions $\mathbb{V}_N$ and $\mathbb{C}_M$ in~\eqref{eq:Kalpha} are those for the samples in~\eqref{eq:vsample}, that is they are given by the compositions
\begin{align}
	\label{eq:eqfcomp}
\mathbb{V}_N(t) := \left( F^{(-1)} \circ \mathbb{U}^V_N \right) (t)
 \quad \mathrm{and} \quad
\mathbb{C}_M(t) := \left( G^{(-1)} \circ \mathbb{U}^C_M \right) (t)  .
\end{align}
So to analyse the asymptotic distribution of   $K_\alpha$,  we will now turn to  the asymptotic behaviour of $\mathbb{E}_\alpha$ specified in~\eqref{eq:Kalpha}
where the terms on the right-hand side are given by~\eqref{eq:eqfcomp}.

Using \eqref{eq:eqfcomp}, \eqref{eq:array} and  conditions~{\bf(A1)}, {\bf(A3)} to take Taylor series expansions with two terms and Lagrange's form of the remainder gives, after elementary transformations,
\begin{align}
	\notag
\mathbb{V}_N (1-t) & = F^{(-1)}\left( 1-t + N^{-1/2} \mathbb{B}^V(1-t) + \theta^V_N(1-t) N^{-1} \ln{N} \right) \\
	\label{eq:taylorF} & = F^{(-1)}(1-t) + \frac{N^{-1/2} \mathbb{B}^V(1-t) + \vartheta^V_N(t) N^{-1} \ln{N}}{f(F^{(-1)}(1-t))}, \quad t \in (0,1),
\end{align}
and, for $  t \in (0,1 \wedge \lambda_\alpha),$
\begin{align}
	\notag
   \mathbb{C}_M (\lambda^{-1}_\alpha t) & = G^{(-1)}\left( \lambda^{-1}_\alpha t + M^{-1/2} \mathbb{B}^C(\lambda^{-1}_\alpha t) + \theta^C_M(\lambda^{-1}_\alpha t) M^{-1} \ln{M} \right)
    \\
	\label{eq:taylorG} & = G^{(-1)}(\lambda^{-1}_\alpha t) + \frac{M^{-1/2} \mathbb{B}^C(\lambda^{-1}_\alpha t) + \vartheta^C_M(t) M^{-1} \ln{M}}{g(G^{(-1)}(\lambda^{-1}_\alpha t))},
\end{align}
where $\vnorm{\vartheta_M^C}_{\infty} + \vnorm{\vartheta_N^V}_{\infty} = O(1)$ a.s.

Substituting  \eqref{eq:taylorF} and \eqref{eq:taylorG}  into the representation for $\mathbb{E}_\alpha$ in \eqref{eq:Kalpha} and using \eqref{eq:ENdef}  gives
\begin{align}
	\label{eq:bbapprox}
 \mathbb{E}_\alpha(t) = E_{\alpha}(t) + N^{-1/2} \mathbb{Z}_{\alpha}(t) + \varphi_{\alpha} (t) N^{-1} \ln N,
\end{align}
where $\vnorm{\varphi_\alpha} = O(1)$ a.s.\ and
\begin{align}
	\label{eq:wnprocess}
\mathbb{Z}_{\alpha}(t)
  :=
  \frac{\mathbb{B}^V(1-t)}{f\left( F^{(-1)} \left( 1 - t \right) \right)}
  -
  \frac{\widehat{\mathbb{B}^C}(\lambda^{-1}_\alpha t)}{ \lambda_\alpha ^{1/2} g\left(\widehat{G^{(-1)}}\left(\lambda^{-1}_\alpha t \right)\right)},
   \quad t \in (0, 1).
\end{align}

One can see from~\eqref{eq:Kalpha1} and \eqref{eq:bbapprox} that if the functional $\mathcal{H}$ were differentiable in a suitable sense at the ``point'' $E_\alpha$, one could derive the desired asymptotic normality of $K_\alpha$ using a suitable version of the functional delta method (see, for example, Ch.~20 of van der Vaart~\cite{vandervaart98}). So we will now turn to studying the local properties of~$\mathcal{H}$ at~$E_\alpha$.

Since $f$, $g$ and $\lambda_\alpha$ are bounded, we see from~\eqref{eq:ENderiv} that 
\begin{align}
\notag -\gamma := \sup_{\alpha\in \mathscr{A} }\sup_{t \in (0,\lambda_\alpha \wedge 1)} E_{\alpha}'(t) <   0  .
\end{align}
As $E_\alpha (t_\alpha)=0,$ the above implies that
\[
E_\alpha (t) \ge  \gamma (t_\alpha - t), \quad 0\le t < t_\alpha ;
\quad
E_\alpha (t) \le \gamma (t_\alpha - t), \quad  t_\alpha < t \le \lambda_\alpha \wedge 1.
\]	
Hence, for any   $v \in \mathscr{D}[0,1],$ setting $\eta := \vnorm{v}_{\infty}$, one has
\begin{align*}
	  \inf_{0\le t < t_\alpha - \eta/\gamma}   (E_\alpha(t_\alpha) + v(t)  )
  & > \gamma \eta/\gamma - \eta = 0,
  \\
	  \sup_{  t_\alpha +  \eta/\gamma <t\le \lambda_\alpha \wedge 1} \left( E_\alpha(t) + v(t) \right)
  & < - \gamma\eta /\gamma + \eta = 0,
\end{align*}
so that the plot of $E_\alpha(t_\alpha) + v(t)$ must ``dive" under zero within $\eta/\gamma$ of $t_\alpha$:
\begin{align*}
	  \left| \mathcal{H}(E_\alpha + v) - \mathcal{H}(E_\alpha) \right|
 \equiv
 \left| \mathcal{H}(E_\alpha + v) - t_\alpha  \right|
 \leq
  \eta/\gamma =   \vnorm{v}_{\infty}/\gamma.
\end{align*}
Thus, the functional $\mathcal{H}$ is Lipschitz continuous in the uniform topology at the   point~$E_{\alpha}$. Using this, \eqref{eq:Kalpha1} and \eqref{eq:bbapprox}, we have
\begin{align}
	\label{eq:deltandef}
 \delta_{\alpha} = O(N^{-1/2})\quad \mbox{a.s.}
\end{align}

\pgfplotsset{
    standard/.style={
        xmin=0,xmax=1.3,
        ymin=-0.6,ymax=0.6,
        axis x line=middle,
        axis y line=middle,
        every axis x label/.style={at={(1,0.5)},anchor=north},
        every axis y label/.style={at={(0,1)},anchor=east}
    }
}

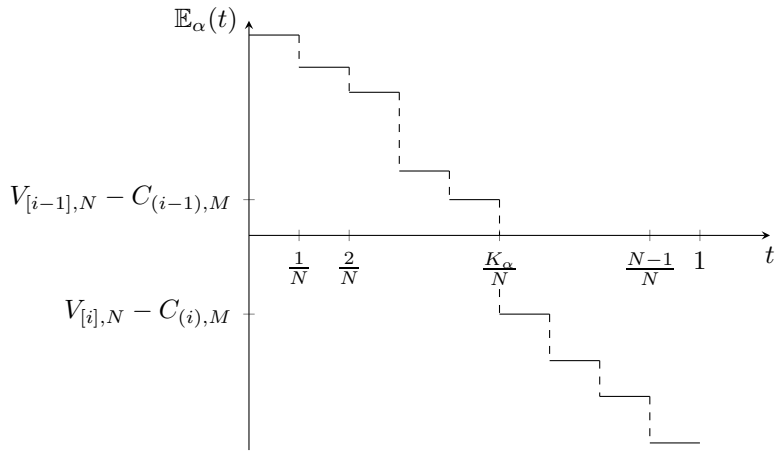
\begin{figure}
	\centering
	\begin{tikzpicture}
	\begin{axis}[
    standard,
    ylabel=\footnotesize{$\mathbb{E}_\alpha(t)$},
    xlabel=\footnotesize{$t$},
    xtick={0,0.125,0.25,0.625,1,1.125},
    xticklabels={\footnotesize{$0$},\footnotesize{$\frac{1}{N}$},\footnotesize{$\frac{2}{N}$},\footnotesize{$\frac{K_\alpha}{N}$},\footnotesize{$\frac{N-1}{N}$},\footnotesize{$1$}},
    ytick={0.1,-0.22},
    yticklabels={\footnotesize{$V_{[i-1],N}-C_{(i-1),M}$},\footnotesize{$V_{[i],N}-C_{(i),M}$}}
]
		\addplot[black, no markers] coordinates {(0,0.56) (0.125,0.56)};
		\addplot[black, dashed, no markers] coordinates {(0.125,0.56) (0.125,0.47)};
		\addplot[black, no markers] coordinates {(0.125,0.47) (0.25,0.47)};
		\addplot[black, dashed, no markers] coordinates {(0.25,0.47) (0.25,0.40)};
		\addplot[black, no markers] coordinates {(0.25,0.40) (0.375,0.40)};
		\addplot[black, dashed, no markers] coordinates {(0.375,0.40) (0.375,0.18)};
		\addplot[black, no markers] coordinates {(0.375,0.18) (0.5,0.18)};
		\addplot[black, dashed, no markers] coordinates {(0.5,0.18) (0.5,0.1)};
		\addplot[black, no markers] coordinates {(0.5,0.1) (0.625,0.1)};
		\addplot[black, dashed, no markers] coordinates {(0.625,0.1) (0.625,0)};
		\addplot[black, dashed, no markers] coordinates {(0.625,-0.15) (0.625,-0.22)};
		\addplot[black, no markers] coordinates {(0.625,-0.22) (0.75,-0.22)};
		\addplot[black, dashed, no markers] coordinates {(0.75,-0.22) (0.75,-0.35)};
		\addplot[black, no markers] coordinates {(0.75,-0.35) (0.875,-0.35)};
		\addplot[black, dashed, no markers] coordinates {(0.875,-0.35) (0.875,-0.45)};
		\addplot[black, no markers] coordinates {(0.875,-0.45) (1,-0.45)};
		\addplot[black, dashed, no markers] coordinates {(1,-0.45) (1,-0.58)};
		\addplot[black, no markers] coordinates {(1,-0.58) (1.125,-0.58)};
\end{axis}
\end{tikzpicture}	
\caption{A graphical illustration of inequalities~\eqref{eq:Eboldr},   \eqref{eq:Eboldl}. Recall that  $K_\alpha N^{-1}= \mathcal{H}(\mathbb{E}_\alpha(t)) = t_\alpha+ \delta_\alpha  $.}
\label{proof}
\end{figure}

Further, since $\mathbb{E}_\alpha$ is right continuous, it follows from the definition of $\mathcal{H}$ that
\begin{align}
	\label{eq:Eboldr} \mathbb{E}_\alpha(t_\alpha + \delta_\alpha) \leq 0.
\end{align}
On the other hand, using the order statistics $V_{[i-1],N} $ (descending) and $C_{(i),M}$ (ascending) for our samples~\eqref{eq:vsample} (see the definitions just before~\eqref{eq:kdef}), we obtain (see Figure~\ref{proof}) for an illustration of the first inequality)
\begin{align}
	\notag \mathbb{E}_\alpha(t_\alpha + \delta_\alpha) & \geq -\max_{2 \leq i \leq K_\alpha} \left[ \left(V_{[i-1],N} - V_{[i],N}  \right) + \left( C_{(i),M} - C_{(i-1),M} \right) \right] \\
	\label{eq:Eboldl} & \geq - \max_{2 \leq i \leq N} \left(V_{[i-i],N} - V_{[i],N} \right) - \max_{2 \leq i \leq M} \left(C_{(i),M} - C_{(i-1),M} \right).
\end{align}
From a well-known result regarding maximal uniform spacings (see, for example, Devroye~\cite{devroye82}) and representation~\eqref{eq:vsample}, we have
\begin{align*}
 \max_{2 \leq i \leq N} \left( V_{[i-1],N} - V_{[i],N} \right) \leq \frac{\max_{2 \leq i \leq N}\left(U_{[i-1],N}^V - U_{[i],N}^V \right)}{\min_{t \in [a,b]} f(t)} = O( N^{-1} \ln{N}) \quad\mathrm{a.s.}
\end{align*}
Since the second term on the right-hand side of~\eqref{eq:Eboldl} has the same order of magnitude,   \eqref{eq:Eboldr} and \eqref{eq:Eboldl}  now yield
\begin{align*}
\mathbb{E}_\alpha(t_\alpha + \delta_\alpha) = O(N^{-1} \ln{N})\quad \mathrm{a.s.}
\end{align*}
Therefore \eqref{eq:bbapprox} implies that
\begin{align}
	\label{eq:Ealphapoint} E_\alpha(t_\alpha + \delta_\alpha) + N^{-1/2}\mathbb{Z}_\alpha(t_\alpha + \delta_\alpha) = O ( N^{-1} \ln{N} ).
\end{align}
Under   assumption {\bf(A3)} we also have
\begin{align}
	\label{eq:tayloren} E_\alpha(t_{\alpha} + \delta_{\alpha}) = E_\alpha(t_{\alpha}) + \delta_{\alpha} E'_\alpha(t_{\alpha}) + O(\delta_{\alpha}^2) = \delta_{\alpha} E'_\alpha(t_{\alpha}) + O(\delta_{\alpha}^2)\quad \mathrm{a.s.}
\end{align}
Combining \eqref{eq:Ealphapoint} with \eqref{eq:tayloren} and using~\eqref{eq:deltandef} gives
\begin{align}
\notag	
  \delta_\alpha  E_\alpha'(t_\alpha)
  & = - N^{-1/2} \mathbb{Z}_\alpha(t_\alpha + \delta_\alpha) + O(N^{-1} \ln{N})
\\
	\label{eq:deltaapp} & = - N^{-1/2} \mathbb{Z}_\alpha(t_\alpha)
 + N^{-1/2}\psi_\alpha \omega_{\mathbb{Z}_\alpha} (\delta_\alpha)  + O(N^{-1} \ln{N})
\end{align}
with $|\psi_\alpha|\le 1$, where
\begin{align}
	\notag \omega_h(\delta) = \sup_{|t-s|\leq \delta} |h(t) - h(s)|
\end{align}
denotes the modulus of continuity of the continuous function~$h$ on $[0,1]$. Recall that, for a Brownian bridge process $\mathbb{B},$ one has
\begin{align*}
\limsup\limits_{\delta \downarrow 0} \frac{w_{\mathbb{B}}(\delta)}{\sqrt{2 \delta \ln{(1/\delta)}}} = 1\quad \mathrm{a.s.}
\end{align*}	
(see, for example, Theorem~1.4.1 in Cs\"{o}rg\H{o} and R\'{e}v\'{e}sz~\cite{cr81}). As this holds for both   $\mathbb{B}^V$ and $\mathbb{B}^C$, and $\lambda_\alpha$, $f$ and $g$ are bounded away from zero, it follows from~\eqref{eq:wnprocess} that
\begin{align}
	\notag \omega_{\mathbb{Z}_\alpha}(\delta_\alpha) = O( \delta_\alpha^{1/2} \ln^{ 1/2}(1/\delta_\alpha))
 = O\bigl(N^{-1/4}(\ln N)^{ 1/2}\bigr) \quad \mathrm{a.s.}
\end{align}
Hence \eqref{eq:deltaapp} now yields
\begin{align}
	\label{eq:deltanfinal}
 \delta_{\alpha} = -\frac{\mathbb{Z}_\alpha(t_{\alpha})}{N^{1/2} E'_\alpha(t_{\alpha})} + O\bigl(N^{-3/4}(\ln N)^{ 1/2}\bigr) \quad \mathrm{a.s.}
\end{align}
\pu{Recall that, for a Brownian bridge process $\mathbb{B}$ and any $t \in (0,1)$, $\mathbb{B}(t)$ is normally distributed with respective mean and variance $(0,t(1-t))$.} Since the processes  $\mathbb{B}^V$ and $\mathbb{B}^C$ in \eqref{eq:wnprocess}  are independent, we immediately see that  $\mathbb{Z}_\alpha(t_{\alpha})$ is normally distributed with zero mean and variance
\begin{align*}
 \frac{t_{\alpha}(1-t_{\alpha})}{f^2(F^{(-1)}(1-t_{\alpha}))} + \frac{t_{\alpha}(1-\lambda^{-1}_\alpha t_{\alpha})}{\lambda^{2}_\alpha g^2(G^{(-1)}(\lambda^{-1}_\alpha t_{\alpha}))}.
\end{align*}	
\pu{Setting}
\[
Z_\alpha^{(1)}:= -  \mathbb{Z}_\alpha(t_{\alpha})/  E'_\alpha(t_{\alpha} ),
\]
\pu{it follows from \eqref{eq:sigmadef} that $Z_\alpha^{(1)}$ is normally distributed with respective mean and variance $(0,\sigma_\alpha^2)$. Thus, \eqref{eq:deltanfinal} establishes the first assertion of Theorem~\ref{thm:1}.}

For welfare, from \eqref{eq:welfQ}, \eqref{eq:Kalpha1}  and \eqref{eq:Kalpha} we have
\begin{align}
	\notag
 W_\alpha
    & =
    \sum_{i=1}^{K_\alpha} \left( V_{[i],N} - C_{(i),M} \right)
    \\
    \notag
    & =
    N \int_{0}^{K_\alpha/N} \left(\mathbb{V}_{N}(1-t) - \mathbb{C}_M(\lambda^{-1}_\alpha t) \right)  dt
    =N \int_{0}^{K_\alpha/N} \mathbb{E}_\alpha (t)\, dt
    \\
     & =
     N \int_{0}^{t_\alpha+\delta_\alpha}
     \bigl(E_{\alpha}(t) + N^{-1/2} \mathbb{Z}_{\alpha}(t) + \varphi_{\alpha} (t) N^{-1} \ln N\bigr) dt.
     \label{eq:w*def}
\end{align}

Now note that   replacing $\int_{0}^{t_\alpha+\delta_\alpha}$ with $\int_{0}^{t_\alpha}$ in the last line   will only introduce an error $O(1)$ a.s. Indeed, since $E_\alpha (t_\alpha) =0$ and   $E_\alpha'(t)$ is bounded in view of our assumptions {\bf (A1)} and {\bf (A2)} (cf.~\eqref{eq:ENderiv}), one has from~\eqref{eq:deltandef} that
\[
N\int_{t_\alpha}^{t_\alpha+\delta_\alpha} E_{\alpha}(t) \, dt
 = O \left(N \int_{t_\alpha}^{t_\alpha+|\delta_\alpha|}
   (t-t_\alpha) \, dt\right)
   = O \left(N\delta_\alpha^2\right) = O(1)\quad \mbox{a.s.}
\]
Further, it is clear from {\bf (A1)}, {\bf (A2)} and~\eqref{eq:wnprocess} that there exists a constant $c<\infty$ such that
\[
 \max_{t\in [0,1]} |\mathbb{Z_\alpha}(t)|\le Y:= c \left( \max_{t\in [0,1]} |\mathbb{B}^V (t)| + \max_{t\in [0,1]}|\mathbb{B}^C (t)|\right)<\infty \quad \mbox{a.s.}
\]
Hence, again using~\eqref{eq:deltandef},
\[
\left|
N^{1/2}\int_{t_\alpha}^{t_\alpha+\delta_\alpha}  \mathbb{Z_\alpha}\, dt
\right|
\le N^{1/2} |\delta_\alpha| Y = O(1)\quad\mbox{a.s.}
\]
This proves the desired claim since the contribution of the last term in the integrand in~\eqref{eq:w*def} to $\int_{t_\alpha}^{t_\alpha+\delta_\alpha}$ will be even smaller in magnitude (recall that $\|\varphi_\alpha\|_\infty=O(1)$~a.s.\ and one has~\eqref{eq:deltandef}).

Thus we obtain  from~\eqref{eq:w*def} that
\begin{align}
\label{eq:Wsimple}
W_\alpha =N \int_{0}^{t_\alpha }
      E_{\alpha}(t)\, dt
      + N^{ 1/2} \int_{0}^{t_\alpha } \mathbb{Z}_{\alpha}(t)\, dt + O(\ln N ) \quad\mbox{a.s.}
\end{align}
It is clear from~\eqref{eq:wnprocess}  that $ \mathbb{Z}_{\alpha}$ is  a zero mean Gaussian process, and so the second integral  on the right-hand side of~\eqref{eq:Wsimple}, to be denoted by $Z_\alpha^{(2)}$, is zero mean normal as well. For brevity, set $f^* (t) := f(F^{(-1)}( t))$ and $g^* (t) := g(G^{(-1)}(t))$. Recalling the right relation in~\eqref{eq:t0def} and also that $\mathbb{B}^V$ and $\mathbb{B}^C$ are independent Brownian bridges, we compute the variance of   $Z_\alpha^{(2)}$ as
\begin{align}
	\notag    \mathbf{E} &   \left( \int_0^{t_\alpha}
   \mathbb{Z}_\alpha(t) \, dt \right)^2
      =
    \iint_{[0,t_\alpha]^2}   \!\!\!
     \mathbf{E}\,  \mathbb{Z}_\alpha(s) \mathbb{Z}_\alpha(t)   \, ds\,  dt
   =
   2 \iint_{[0,t_\alpha]^2\cap\{s<t\}} \hspace{-1cm}
    \mathbf{E} \,  \mathbb{Z}_\alpha(s) \mathbb{Z}_\alpha(t)   \, ds\, dt
   \\
	\notag &
  = 2 \int_0^{t_\alpha} \int_0^{t} \left[
   \frac{\mathbf{E}\,  \mathbb{B}^V(1-s) \mathbb{B}^V(1-t)
   }{f^* (1-s) f^*(1-t)}
    +
   \frac{
    \mathbf{E}\, \mathbb{B}^C(\lambda^{-1}_\alpha s) \mathbb{B}^C(\lambda^{-1}_\alpha t)
    }{ \lambda_\alpha  g^*( \lambda^{-1}_\alpha s) g^*( \lambda^{-1}_\alpha t)}
    \right] \, ds \, dt
   \\
	\notag
 & = 2 \int_0^{t_\alpha} \int_0^{t} \left[ \frac{s(1-t)}{f^* (1-s) f^*(1-t)}
 +
 \frac{s(1 - \lambda^{-1}_\alpha t) }{ \lambda_\alpha^2  g^*( \lambda^{-1}_\alpha s) g^*( \lambda^{-1}_\alpha t)} \right] \, ds \, dt
 \\
 \label{eq:Wvar}
 & = 2 \int_0^{t_\alpha}\left[ \frac{1-t}{f^*(1-t)}
  \int_0^{t}
  \frac{s\, ds  }{f^* (1-s) }
 +
 \frac{1 - \lambda^{-1}_\alpha t}{ g^*( \lambda^{-1}_\alpha t)}
 \int_0^{t}
 \frac{s \, ds }{ \lambda_\alpha^2  g^*( \lambda^{-1}_\alpha s)}  \right] dt  .
\end{align}
Changing variables $x:=F^{(-1)}( 1-s)$ and $x:=G^{(-1)}( \lambda^{-1}_\alpha  s)$, respectively, in the two integrals inside the square brackets in the last line in~\eqref{eq:Wvar}, we obtain the desired representation~\eqref{eq:varsigmadef} for~$\varsigma_\alpha^2.$ Now the second assertion of the theorem follows from~\eqref{eq:Wsimple}.

To complete the proof of the theorem, it remains to compute the covariance
\begin{align*}
 \mathbf{E}\, Z_\alpha^{(1)}Z_\alpha^{(2)}
   = -\mathbf{E}\,    \frac{\mathbb{Z}_\alpha(t_{\alpha})}{  E'_\alpha(t_{\alpha} )} \int_{0}^{t_\alpha } \mathbb{Z}_{\alpha}(t)\, dt
    = -  \frac{1}{  E'_\alpha(t_{\alpha} )}
    \int_{0}^{t_\alpha }{\mathbf E}\, \mathbb{Z}_\alpha(t_{\alpha}) \mathbb{Z}_{\alpha}(t)\, dt
    = -  \frac{S_\alpha (t_\alpha)}{  E'_\alpha(t_{\alpha} )},
\end{align*}
where the last equality follows from evaluation of the inner integral in~\eqref{eq:Wvar} and the definition of the function~$S_\alpha$ (following~\eqref{eq:varsigmadef}). Theorem~\ref{thm:1} is proved.
\eop

\begin{remark} \em{An argument similar to the second part of our proof of the asymptotic normality of $K_\alpha$ can be found in the first example in Section~8, Ch.~I of Borovkov~\cite{borovkov98}. However, in our case the quantity $t_{\alpha}$ depends on~$\alpha$ and  does not converge to any value. Therefore, the first-order derivative of $\mathcal{H}$ cannot be computed by directly applying a known general version of the delta method, such as the main result in Borovkov~\cite{borovkov85}.} \end{remark}

\begin{remark} \em{The reader is referred to Cs\"{o}rg\H{o} and R\'{e}v\'{e}sz~\cite{cr78} for the minimal boundedness conditions that must be imposed on $f'$ and $g'$ for \eqref{eq:taylorG} and \eqref{eq:taylorF} to hold. It turns out that, for the approximation of $K_\alpha$, these boundedness conditions may be weakened further as they are required only in a neighbourhood of $t_{\alpha}$.} \end{remark}

\medskip
\noindent
{\em Proof of Theorem~\ref{thm:2}.} (i)~When $G = F,$ the equation $E_\alpha (t_\alpha)=0$ for $t_\alpha$ (see~\eqref{eq:ENdef}, \eqref{eq:t0def}) turns into $ F^{(-1)}(\lambda^{-1}_\alpha t_\alpha ) = F^{(-1)}(1 - t_\alpha )$, which means that $\lambda^{-1}_\alpha t_\alpha  = 1 - t_\alpha $, yielding the desired representation for~$ t_{\alpha}.$

Next, in this case for $t=t_\alpha$ the value of~\eqref{eq:ENderiv} turns into
\[
E_\alpha' (t_\alpha) = -    \frac{1 }{t_\alpha f(F^{(-1)} (1-t_\alpha)) },
\]
so that~\eqref{eq:sigmadef} becomes
\[
\sigma_\alpha^2 = t_\alpha^2
\bigl[ t_\alpha(1 -t_\alpha) + t_\alpha ^2 \lambda_\alpha^{-2}\bigr]
 = t_\alpha^2 (1-t_\alpha^2) = \frac{\lambda_\alpha^2}{(1+\lambda_\alpha)^3}.
\]

(ii)~As we already pointed out, in the special  case when $F=G$ the exact distribution of~$K_\alpha$ does not depend on~$F$. In fact, as shown in Loertscher, Muir and Taylor~\cite{lmt16a}, $K_\alpha$ has the hypergeometric distribution $\mathrm{Hg}(N,M, N + M)$. It is well known that, as $M+N\to \infty$ and the quantity $\lambda_\alpha$ remains bounded away from~0 and~1 (which is ensured by our assumption~{\bf (A2)}), that distribution can be approximated by a normal distribution (see, for example, Theorem~2.1 in Lahiri and Chatterjee~\cite{lc07}). The stated bounds \eqref{eq:BerryEsseen} follow  from the results established in Theorem~2.2  in Lahiri and Chatterjee~\cite{lc07}.

(iii)~First note that, in the case of the uniform distributions $F=G$ on $(0,1)$, one has $F^{(-1)}(t) \equiv G^{(-1)}(t)\equiv t$ on $(0,1)$, so that
\begin{align*}
\int_0^{t_\alpha} E_\alpha(t)\,dt
 = \int_0^{t_\alpha} \left(1 - t- t/\lambda_\alpha\right) \, dt
  =  \int_0^{t_\alpha} \left(1 -   t/t_\alpha\right) \, dt =\frac{t_\alpha}{2}=\frac{\lambda_\alpha}{2(1+\lambda_\alpha)}.
\end{align*}
Next, since $f(t)\equiv g(t) \equiv 1$ on $(0,1)$, we also have
\begin{align*}
S_\alpha (t) &
 = (1-t) \int_{1-t}^1 (1-x)\, dx
+ (1-\lambda_\alpha^{-1}t) \int_{0}^{\lambda_\alpha^{-1}t} x\, dx
 \\
& = \frac12 t^2 (1-t)
 + \frac12 \left(\frac{t}{\lambda_\alpha}\right)^2
 \left(1-\frac{t}{\lambda_\alpha}\right), \quad t\in (0,t_\alpha).
\end{align*}
Integrating this expression in $t$ from 0 to $t_\alpha$ yields the second formula in~\eqref{eq:varsigma}. To get the last formula in~\eqref{eq:varsigma}, we just note that $E_\alpha (t_\alpha)'=-1/t_\alpha$ and, as $ t_\alpha/\lambda_\alpha= 1-t_\alpha,$ one has   $S_\alpha (t_\alpha)= \frac12 t_\alpha^2 (1-t_\alpha)
 + \frac12 (1-t_\alpha)^2 t_\alpha = \frac12 t_\alpha  (1-t_\alpha)$, and so $\varkappa_\alpha=- S_\alpha (t_\alpha)/E_\alpha (t_\alpha)'=\frac12 t_\alpha^2 (1-t_\alpha).$
\eop

\section*{Acknowledgements}

E.~Muir's research  was supported by the Elizabeth and Vernon Puzey Scholarship and the Australian Research Council through the Laureate Fellowship FL130100039 and the ARC Centre of Excellence for Mathematical and Statistical Frontiers. K.~Borovkov's work  was supported by the ARC grant  DP150102758.

\setlength{\baselineskip}{0.67\baselineskip}

\end{document}